\documentclass[12pt]{amsart}

\usepackage{amsmath,amsthm,amssymb,eucal}

\newtheorem{thm}{Theorem}
\newtheorem{cor}[thm]{Corollary}

\theoremstyle{definition}
\newtheorem{defn}{Definition}
\newtheorem{rk}{Remark}
\newtheorem{Claim}{Claim}

\theoremstyle{remark}

\def\dfrac#1#2{{\displaystyle \frac{#1}{#2}}}

 \def\Rset{{\mathbb R}}

\numberwithin{equation}{section}

\begin{document}

\title{
Nonsmooth Critical Point Theorems\\
Without Compactness
}
\author{
		Youssef JABRI
		}
\address{
     University Mohamed~I, Department of Mathematics,     
     Faculty of Sciences, Box 524, 60000 Oujda, Morocco
     }
\email{jabri@sciences.univ-oujda.ac.ma} 
\date{}

\begin{abstract}
We establish an abstract critical point theorem for locally Lipschitz 
functionals that does not require any compactness condition of 
Palais-Smale type.  It generalizes and unifies three other critical 
point theorems established in \cite{J-M} for $C^{1}$-functionals under 
slightly stronger assumptions.
Our approach uses continuous selections of multivalued mappings.
\end{abstract}

\maketitle



\bigskip  

Abstract critical point theorems from the linking family like the 
mountain pass theorem of Ambrosetti-Rabinowitz and saddle point 
theorem of Rabinowitz, (cf.  for example 
\cite{Rabinowitz,M-W-book,Struwe}) are important tools in nonlinear
analysis.  They require generally a compactness condition known as 
Palais-Smale condition, (PS) for short.  In applications to nonlinear 
boundary value problems, verifying if (PS) holds or not is crucial and 
turns out often to be technical and/or too long.

In \cite{J-M}, the author and Moussaoui proved a critical point 
theorem (and some variants) with similar conditions to Rabinowitz 
saddle point theorem without requiring the Palais-Smale condition.  
They supposed instead a convexity assumption on $\Phi$ on a part of 
the underlying space.

In this paper, we are concerned with a twofold extension of the 
results of \cite{J-M}.  We will show that they follow easily from a 
more abstract critical point theorem valid for locally Lipschitz 
functionals with ``less convexity'' on $\Phi$.  The proof also is 
easier.
A novelty in our approach is that we exploit the notion of 
\emph{continuous selections of multivalued mappings}.

The following result is a reliable prototype of our critical 
point theorem without compactness.
\begin{thm}
Let $E$ be a reflexive banach space such that $E=V \oplus W$ where 
$\dim V<+\infty$ and $\Phi\colon E\to \Rset$ a $C^1$-functional that 
satisfies:\\
\textup{(i)} $\Phi$ is coercive on $W$, i.e., $\Phi(w)\to +\infty$ as 
$\Vert w\Vert \to \infty$.\\
\textup{(ii)} $\forall w\in W$, $v\mapsto \Phi(v+w)$ is 
quasi-concave.\\
\textup{(iii)} $\Phi(v+w)\to -\infty$ as $||v||\to +\infty$ and the 
convergence is uniform on bounded subsets of~$W$.\\
\textup{(iv)} $\forall v\in V$, $w\mapsto \Phi(v+w)$ is weakly lower 
semi-continuous.

Then, $\Phi$ admits a critical point $u$ in $E$ characterized by the 
minimax argument
$$
c=\min_{w\in W} \max_{v\in V} \Phi(v+w).
$$\label{thm1}
\end{thm}

In \cite{J-M}, it is proved in the case of a \emph{Hilbert space} $E$, 
under the hypotheses (i), (iii), (iv) and
$$
\leqno{\mbox{(ii$'$)}} \quad \forall w\in W,\ v\mapsto \Phi(v+w) 
\text{ is concave.}
$$
that the minimax value $c$ is a critical value. 
The concavity and the Hilbertian structure of $E$ play a crucial role 
in the proof.

In the case of a reflexive Banach space, it was proved (more easily) 
that the same conclusion holds under the hypotheses (i), (iii), (iv) 
and
$$
\leqno{\mbox{(ii$''$)}} \quad \forall w\in W,\ v\mapsto \Phi(v+w) 
\text{ is strictly quasi-concave.}
$$
In \cite[Remark 15]{J-M}, they pointed out that in the Hilbertian 
case, the question wether the quasi-concavity was sufficient or no to 
guarantee the existence of a critical point was open.
This is intimately related to the fact that with (i), (iii), (iv) and 
(ii$''$), $\Phi$ does not satisfy in general the local Palais-Smale 
condition (PS)$_c$ for the minimax critical value $c$.
While, as they mentioned in \cite[Remark 12]{J-M}, with (ii$'$), it is 
unclear if $\Phi$ satisfies (PS)$_c$ or not.

By Theorem~\ref{thm1}, we know now that the quasi-concavity is 
sufficient to get a critical value, even in the case of a reflexive 
Banach space.

The paper is organizd as follows.  In Section 2, we describe the 
origins of our minimax theorems.  Then, in Section 3, we recall the 
necessary background from critical point theory for locally Lipschitz 
functionals following Clarke~\cite{Clarke}.  We recall also some 
elements of topology related to continuous selections of multivalued 
mappings.  Section 4 is devoted to our main result and its differents 
corollaries generalizig all results of~\cite{J-M}.

\section{The Origins of the Method}

The origins of the minimax procedure that appears in 
Theorem~\ref{thm1} goes back to \textbf{1949}.

Let $\Omega$ be a smooth bounded domain of $\Rset$, $f\colon \Omega 
\times\Rset \to \Rset$ a Carath\'eodory function.  Consider the 
nonlinear Dirichlet problem
$$
(\mathcal{P})\qquad \left\{
\begin{array}{rll}
-\Delta u & =f(x,u) & \text{in }\Omega,\\
u& =0 & \text{on }\partial \Omega.
\end{array}
\right.
$$
and denote the potential associated to $f$ by 
$$
F(x,t)=\int_0^t f(x,s)\,ds.
$$
Dolph \cite{Dolph} solved the problem $(\mathcal{P})$ when
\begin{equation}
\lambda_k < \mu_k \leq \liminf_{s\to \pm \infty} \frac{f(x,s)}{s} 
\leq 
\limsup_{s\to \pm \infty} \frac{f(x,s)}{s} \leq \mu_{k+1} < 
\lambda_{k+1}
\label{eq1}
\end{equation}
where $\lambda_k$ and $\lambda_{k+1}$ are two consecutive eigenvalues 
of $-\Delta$ in $H^1_0(\Omega)$ (nonresonance between two consecutive 
eigenvalues of the Laplacian).

A similar condition in terms of the potential $F$ may be expressed as
\begin{equation}
\lambda_k < \mu_k \leq \liminf_{s\to \pm \infty} \frac{2F(x,s)}{s^2} 
\leq \limsup_{s\to \pm \infty} \frac{2F(x,s)}{s^2} \leq \mu_{k+1} < 
\lambda_{k+1}
\label{eq2}
\end{equation}

In general, variational methods fail to handle the problem 
$(\mathcal{P})$ when only conditions on the potential, like 
\eqref{eq2}, are required (cf.  the discussion of \cite{C-O} for 
example).  This makes it impossible to verify the Palais-Smale 
condition required in linking theorems.

The first variational attempt to solve $(\mathcal{P})$ under condition 
\eqref{eq2}, without requiring any assumption on $f$, is also due to 
Dolph.  He required the following additional condition
$$
\textup{(DO)}\qquad \left\{
\begin{array}{l}
\text{The energy functional associated to the problem }(\mathcal{P})\\
\displaystyle\qquad\qquad \Phi(u)=\frac{1}{2}\int_\Omega|\nabla u|^2 
- 
\int_\Omega F(x,u)\,dx,\\
\text{admits at most one maximum in }w+V\text{ for all }w\in W.
\end{array}
\right.
$$
The spaces $V=\bigoplus_{i\leq k} E(\lambda_{i})$ and $W=V^{\perp}= 
\bigoplus_{i\geq k+1} E(\lambda_{i})$ where $E(\lambda_{i})$ is the 
eigenspece associated associated to the eigenvalue $\lambda_{i}$ of 
the Laplace operator with Dirichlet boundary data.

Later, Thews \cite{thews} treated $(\mathcal{P})$ under condition 
\eqref{eq2} allowing $\mu_k=\lambda_k$ and $\mu_{k+1}=\lambda_{k+1}$ 
but supposed more restrictive conditions than~(DO).

Theorem \ref{thm1} may be seen as a continuation of the attempt of the 
author and Moussaoui in \cite{J-M} to provide an abstract 
framework to Dolph's method.  In \cite{J-M}, this approach was applied 
to get a variant of a result proved earlier by Mawhin and Willem 
in~\cite{M-W-paper} by combining the least dual action of 
Clarke-Ekeland \cite{Clarke-Ekeland} and an approximation method of 
Br\'ezis~\cite{Brezis}.

\section{Some Elements of Nonsmooth Critical Point Theory and 
Multivalued Analysis}

Let $\Phi\colon X\to \Rset$ be a locally Lipschitz functional. For 
each $x$, $v\in X$, the generalized directional derivative 
of $\Phi$ at $x$ in the direction $v$ is 
$$
\Phi^{\circ}(x;v)=\limsup_{y\to x,\: t\downarrow 0}
\frac{\Phi(y+tv)-\Phi(y)}{t}.
$$
It follows by the definition of locally Lipschitz functionals that 
$\Phi^{\circ}(x;v)$
is finite and $|\Phi^{\circ}(x;v)| \leq C||v||$.

Moreover, $v\mapsto \Phi^{\circ}(x;v)$ is positively homogenous and 
subadditive and $(x,v) \mapsto \Phi^{\circ}(x;v)$ is u.s.c.

The generalized gradient (Clarke subdifferential) of $\Phi$ at $x$ is 
the subset $\partial \Phi(x)$ of $X^*$ defined by
$$
\partial \Phi(x)=\big\{x^*\in X^*;\ \Phi^{\circ}(x;v)\geq \langle 
x^*, v\rangle,\text{ for all }v\in X \big\}.
$$
It enjoys the following properties:\\
a) For each $x\in X$, $\partial \Phi(x)$ is non-empty, convex weak-$*$ 
compact subset of $X^*$.\\
b) For each $x$, $v\in X$, we have 
$$
\Phi^{\circ}(x;v)=\max\{\langle x^*, v\rangle;\ x^*\in \partial 
\Phi(x) \}
$$ 
c) $\partial (\Phi+\Psi)(x) \subset \partial (\Phi) + \partial 
(\Psi)$, where $\Phi$ and $\Psi$ are locally Lipschitz at~$x$.
\begin{thm}[Lebourg mean value theorem]
If $x$ and $y$ are two distinct points in $X$, then there exists 
$z=x+\tau (y-x)$, $0<\tau<1$ such that
$$
\Phi(y) - \Phi(x) \in \langle \partial \Phi(z), y-x\rangle.
$$
\end{thm}
The notion of critical point of a locally Lipschitz functional is the 
following.
\begin{defn}
Let $\Phi\colon X\to\Rset$ be locally Lipschitz.  A point $x\in X$ is 
a critical point of $\Phi$ if $0\in \partial \Phi(x)$.  A real number 
$c$ is called a critical value of $\Phi$ if $\Phi^{-1}(c)$ contains a 
critical point $x$.
\end{defn}

We recall now some results on multivalued mapppings.
Let $M$ and $N$ be two topological spaces. 
\begin{defn}
A multivalued mapping $T\colon M\to 2^N$ is a map which assigns to 
each point $m\in M$ a subset $T(x)$ in $N$.

A multivalued mapping $T\colon M\to 2^N$ is upper semi-continuous if 
and only if $T^{-1}(A)$ is closed for all closed subsets $A$ of $N$, 
where the preimage $T^{-1}(A)$ is defined by
$$
T^{-1}(A)=\big\{m\in M;\ T(m)\cap A\not= \varnothing \big\}.
$$  
The multivalued mapping $T$ is lower semi-continuous if and only if 
$T^{-1}(A)$ is open for all open subsets $A$ of $N$.

And $T$ is contiuous if and only if it is both lower and upper 
semi-contiunuous.
\end{defn}

\begin{rk}
Notice that when $T$ is single-valued, the lower semi-continuity 
(resp.  upper semi-continuity) of $T$ as multivalued coincids with 
continuity.
\end{rk}

Let $T\colon M\to 2^N$ be a multivalued mapping.  By a 
\emph{selection} of $T$, we mean a single-valued mapping $s\colon M\to 
N$ with
$$
s(m)\in T(m),\qquad \text{ for all }m\in M.
$$
\begin{thm}[Michael's selection theorem]
A lower semi-continuous multivalued mapping $T\colon M\to 2^N$ has a 
\emph{contiunous selection} $s\colon M\to N$ if the following three 
conditions are satisfied:\\
\textup{(i)} $M$ is paracompact,\\
\textup{(ii)} $N$ is a Banach space,\\
\textup{(iii)} The set $T(m)$ is nonempty, closed and convex for all 
$m\in M$.
\end{thm}
\begin{thm}[Minimal selection theorem]
Let $T\colon M\to 2^N$ be a continuous multivalued mapping, where $M$ 
is a metric space and $N$ is a Hilbert space.  Suppose that $T(m)$ is 
nonempty, closed and convex for all $m\in M$.  Denote by 
$\mathfrak{m}(T(m))$ the unique element of the set $T(m)$ with 
smallest norm.

Then $\mathfrak{m}\colon M\to N$ is a \emph{contiunous selection}.
\end{thm}

\section{Continuous Selections of Multivalued Mappings}

We give now our main abstract critical point theorem that contains all the 
subsequent results as special cases.
\begin{thm}
	Let $E$ be a Banach space such that $E=V\oplus W$ and $\Phi\colon 
	E\to \Rset$ a locally Lipschitz functional. Suppose the 
	following assumptions.\\
	\textup{(a)} $\forall w\in W$, the set 
	$$
	V(w)=\big\{v\in V;\ 
	\varphi(w)=\Phi(v+w)=\max_{g\in V}\Phi(g+w) \big\} \not= 
	\varnothing.
	$$
    \textup{(b)} The functional $\varphi\colon W\to\Rset$ is bounded below 
    and achieves its minimum at some point~$\overline{w}$.\\
    \textup{(c)} There exists a continuous selection $s\colon W\to V$ such 
    that $s(w)\in V(w)$ for all $w\in W$.
	
    Then, $\overline{u}= s(\overline{w}) + \overline{w}$ is a critical 
    point such that
	$$
	\Phi(\overline{u}) = \min_{w\in W} \max_{v\in V} \Phi(v+w).
	$$
	\label{thm-multi}
\end{thm}
\begin{proof}
	Take $g\in V$, then
	$$
    \Phi(\overline{u}+t(-g)) - \Phi(\overline{u}) \leq 0,\qquad \forall 
    t>0.
	$$
	Divide by $t$ and let $t$ go to infinity to obtain
	$$
	\Phi^\circ(\overline{u};-g)  = -\Phi^\circ(\overline{u};-g) \leq 0.
	$$
	This is true for all $g$ in the linear space $V$, so
	$$
	\Phi^\circ(\overline{u};g)  = 0,\qquad \text{for all }g\in V.
	$$

    On the other side, we know by Lebourg mean value theorem, that
    $$
    \Phi(\overline{u}+sk) - \Phi(\overline{u}) = \langle \xi,sk\rangle = 
    s\langle \xi,k\rangle \leq s \Phi^\circ(z;k)
    $$
    where $\xi\in \partial \Phi(z)$ for some $z=\overline{u}+\tau w$, 
    $0<\tau<1$ and $k\in X$.
        
    So, for $h\in W$, if we write $w_{t} = \overline{w} + th$, $0<t\leq 1$ 
    and $v_{t}= s(w_{t})$.
	
    Consider a sequence $t_{n}\downarrow 0$ and denote by $\overline{v} = 
    s(\overline{w}) = \lim_{n\to\infty} s(w_{t_{n}})$.\\
    Since $\Phi(\overline{w}+ \overline{v}) \geq \Phi(\overline{w}+ 
    s(w_{t_{n}}))$ because $\overline{v}\in V(\overline{w})$, we have
	$$
    \dfrac{\Phi(\overbrace{\overline{w}+t_{n}h}^{w_{t_{n}}}+ v_{t_{n}}) 
    -\Phi(\overline{w}+v_{t_{n}})}{t_{n}} \geq 
    \dfrac{\Phi(\overbrace{\overline{w}+t_{n}h}^{w_{t_{n}}}+ v_{t_{n}}) 
    -\overbrace{\Phi(\overline{w}+\overline{v})}^{\inf \varphi}}{t_{n}} 
    \geq 0.
	$$
	So, 
	$$
    \Phi^\circ(z_{n},h)\geq 0,\quad\text{where }z_{n} \in 
    ]\overline{w}+v_{t_n}, \overline{w}+v_{t_n}+t_n h[.
	$$
	At the limit we get by the u.s.c. of $\Phi^\circ(.,.)$,
	$$
    \Phi^\circ(\overline{u},h)\geq\limsup_{n\to\infty}\Phi^\circ(z_{n},h) 
    \geq 0.
	$$
	And there too, we get 
	$$
	\Phi^\circ(\overline{u};h)  = 0, \qquad \text{for all }h\in W.
	$$
    So, we have finally 
    that
    $$
    0\in \partial \Phi(\overline{u}).
    $$
\end{proof}
Using this result and Michael's selection theorem, we obtain the 
following immediate consequence.
\begin{cor}
	Suppose that $E$ is as in the former theorem and that $\Phi$ is 
	locally Lipschitz on $X$ and satisfies (a) and (b).
	
    Suppose also that the multivalued mapping, $T\colon \to 2^{V}$, 
    $w\mapsto V(w)$ is lower semi-continuous and $T(w)$ is convex for all 
    $w\in W$.  Then, the conclusion of Theorem \ref{thm-multi} holds true.
\end{cor}
\begin{cor}
	Suppose that $E$ and $\Phi$ are as above and satisfy \textup{(i)}, 
	\textup{(ii$''$)}, \textup{(iii)} and \textup{(iv)}.
	
    Then, the conclusion of Theorem \ref{thm-multi} holds true.
\end{cor}
This result has been proved in a direct way in \cite[page 372]{J-M}.
\begin{proof}
    To see that it is a consequence of Theorem~\ref{thm-multi}, it 
    suffices to show that the single valued $w\mapsto s(w)$ where $s(w)$ 
    is the unique element in $V$ (by strict quasi-concavity) that achieves 
    that maximum of $v\mapsto \Phi(v+w)$ in~$V$ is continuous.
    Consider a sequence $w_{n}\to \overline{w}\in W$.  
    There exists $A>0$, such that $||v||\geq A$ implies (by (iv)) that
    $$
    \Phi(v+\overline{w}) \leq \Phi(v+\overline{w}) < \inf_W \Phi\leq \Phi(w_t).
    $$
    So, $||s (w_t)|| \leq A$.  Because otherwise we get the 
    following contradiction with the definition $v_\varepsilon (w_t)$:
    $$
    \Phi(\varphi_\varepsilon (w_t) + w_t) < \Phi(w_t).
    $$
    Hence, there is a sequence $t_n \to 0$ such that $\varphi_\varepsilon 
    (w_{t_n}) \to v_{0} \in V$.

    While by definition of $s(w_n)$, we have
	$$
	\Phi(s(w_{{n}}) + w_{{n}})\geq\Phi(v+w_{{n}}),\ \forall v\in V.
	$$
	At the limit, we get
	$$
	\Phi(v_0+ w_0) \geq \Phi(v+w_0),\ \forall v\in V,
	$$
	that is, $v_0 =s (w_0)$.
\end{proof}
\begin{cor}
	Suppose that $E$ is a Hilbert space and $\Phi$  satisfies 
	\textup{(i)}, \textup{(ii$"$)}, \textup{(iii)} and \textup{(iv)}.
	
    Then, the conclusion of Theorem \ref{thm-multi} holds true.
\end{cor}
There too, it suffices to show that there is a continuous selection of 
$T$.  It may come to mind to use the minimal selection theorem.  But 
unfortnately, if the sets $V(w)$ are convex and compact (not only 
closed), $T$ is only upper semi-continuous.  Our hypotheses do not 
suffice to get that $T$ is lower semi-continuous.
Counter-examples exist.  Nevertheless, we are able to use a 
perturbation idea seen above to prove the minimal selection is indeed 
continuous in our case.
\begin{proof}
    Consider the single valued mapping $w\in W\mapsto s(w)= \max_{V(w)} 
    (\Phi(v+w) - ||v||)$.  It is well defined because, the sets $V(w)$ are 
    convex (by (ii$'$) and closed (they are even compact by (iii) and 
    $\dim V<+\infty$).  And $s(w)$ is indeed the minimal selection 
    because in $V(w)$, $\Phi(v+w)$ is constant and the element 
    with smallest norm of the closed convex set $V(w)$ in the Hilbert 
    $E$, 
    is unique.  Moreover it realizes the maximum of $v\in V(w)\mapsto 
    \Phi(v+w) - ||v||$.
	
	Let us show that $s$ is continuous in $W$. Consider $w_{n}\to 
	\overline{w}$. Then, $s(w_{n})$ is bounded. Otherwise, there would 
	be $N>0$ such that  for all $n\geq N$,
	$$
    \Phi(s(w_{n})+w_{n}) - ||s(w_{n})|| < \inf_W \Phi\leq \Phi(0+w_n)- 
    ||0||=\Phi(w_{n}).
    $$
    A contradiction since $s(w_{n})$ achieves the maximum of $\Phi$ 
    on~$w_{n}+V$.
	So, $s(w_{n})\to \overline{v}$.
	
	First $\overline{v}\in V(\overline{w})$, indeed
	$$
	\begin{array}{rl}
		\Phi(\overline{v}+\overline{w}) & =\lim_{n} \Phi(w_{n}+s(w_{n}))  \\
		 & \geq \lim_{n} \Phi(w_{n}+v),\qquad \forall v\in V  \\
		 & \geq  \Phi(\overline{w} +v),\qquad \forall v\in V.
	\end{array}
	$$
	On the other hand, we have  by definition $s(w_n)$, 
	$$
    \Phi(s(w_{{n}}) + w_{{n}}) - ||s(w_{n})||\geq\Phi(v+w_{{n}}) - ||v||,\ 
    \forall v\in V.
	$$
	And at the limit, we get
	$$
	\Phi(w_0+ v_0) - ||v_{0}||\geq \Phi(v+w_0),\ \forall v\in V,
	$$
	and $v_0 =s (w_0)$.
\end{proof}

Finally, we come to the proof of Theorem \ref{thm1} which is in fact 
also a simple corollary of Theorem~\ref{thm-multi}.
\begin{proof}[Proof of Theorem \ref{thm1}]
    The proof of the former corollary applies verbatim except for the 
    justification of the definiteness of ``the single valued'' 
    mapping~$s$.  In the proof of the former corollary (Hilbertian case), 
    we used the fact that $V(w)$ is convex and closed.  So, a unique 
    element of minimal norm exists.
    In our new situation, this follows from the strict quasi-convexity of 
    $v\mapsto \Phi(v+w) - \Vert v\Vert$ (by (ii) and the strict convexity 
    of the norm in a reflexive Banach space \cite{asplund}) and from the 
    \emph{compactness} of $V(w)$ for each $w\in W$.
\end{proof}

We can remove the assumption that $\dim V<+\infty$ by requiring that 
$\Phi$ is weakly upper semi-continuous.  This improves Theorem~11 in 
\cite{J-M} and shows that it is also a particular case of 
Theorem~\ref{thm-multi}.
\begin{thm}
    In a Banach reflexive space $X=V\oplus W$ such that (i), (ii), (iii) 
    are satisfied.  Suppose moreover that $\Phi$ is weakly upper 
    semi-continuous.  Then, the conclusion of Theorem~\ref{thm1} holds 
    true.
	\label{thm-infinite}
\end{thm}
\begin{proof}[Proof of Theorem \ref{thm-infinite}]
    The reflexive character of $V$ with the anti-coerciveness of $v 
    \mapsto \Phi(v+w)$ suffice to guarantee that $V(w)$ is nonempty and 
    bounded.  Since it is also closed and convex (by (ii)) it is weakly 
    compact.  So, the strict quasi-concave weakly upper semi-continuous 
    functional $v \mapsto \Phi(v+w) -\Vert v\Vert$ achieves its maximum in 
    a unique point $s(w)$.
	\begin{Claim}
		The selection $w\in W \mapsto s(w)\in V(w)$ is continuous.
	\end{Claim}
    Consider a sequence $(w_{n})_{n}\subset W$ such that $w_{n}\to 
    \overline{w}$.  Then, the sequence $v_{n}=s(w_{n})$ is bounded.
	Indeed, by (iv), there exists $A>0$ such that
	$$
	\Phi(v+w_{n}) - \Vert v\Vert < \inf_{W} \Phi < \Phi(w_{n}) = 
	\Phi(w_{n}+ 0) - \Vert 0\Vert,\quad \forall v\in V,\ \Vert v\Vert 
	\geq A.
	$$
	So, $\Vert v_{n} \Vert \leq A$, because otherwise we would get the 
	contradiction
	$$
	\Phi(v_{n}+w_{n}) < \Phi(w_{n}).
	$$
	Therefore, $v_{n}\rightharpoonup \overline{v}\in V$.
	This weak limit $\overline{v}$ belongs to $V(\overline{w})$. Indeed,
	$$
	\begin{array}{rl}
		\Phi(\overline{v}+\overline{w}) & \geq 
                 \limsup_{n}\Phi(v_{n}+w_{n})  \\
		 & \geq \limsup_{n}\Phi(v+w_{n}),\qquad \forall v\in V  \\
		 & \geq \Phi(v+\overline{w}),\qquad \forall v\in V.
	\end{array}
	$$
	Moreover, we know that
	$$
    \Phi(v_{n}+w_{n}) - \Vert v_{n}\Vert\geq \Phi(v+w_{n}) - \Vert 
    v\Vert,\qquad \forall v\in V.
	$$
    By the weak upper semi-continuity of $\Phi$, the weak lower 
    semi-continuity of the norm and since $w_{n}+ v_{n} \rightharpoonup 
    \overline{w} + \overline{v}$,
	$$
    \Phi(\overline{v}+\overline{w})- \Vert \overline{v}\Vert \geq 
    \Phi(v+\overline{w}) - \Vert v\Vert,\qquad \forall v\in V,
	$$
	i.e., $\overline{v}=s(\overline{w})$.
\end{proof}
Some applications to nonlinear boundary value problems will be 
investigated in a forthcoming paper.

\end{document}